\documentclass[11pt,twoside]{article}
\title{Deterministic particle method approximation of a 
contact inhibition cross-diffusion problem
\thanks{First author supported by the Spanish MCINN Project
MTM2010-18427. Second author supported by the Spanish MCINN Project MTM2010-21135-C02-01.}}

\author{Gonzalo Galiano  \thanks{Dpt. of Mathematics, Universidad de Oviedo,
 c/ Calvo Sotelo, 33007-Oviedo, Spain ({\tt galiano@uniovi.es, selgasvirginia@uniovi.es})}
    \and Virginia Selgas\footnotemark[2] }
\date{}
\usepackage{fancyhdr}
\usepackage{amsfonts}
\usepackage{amsmath}
\usepackage{amssymb}
\usepackage{amsthm}
\usepackage{graphics} 
\usepackage{epsfig} 
\usepackage{latexsym}
\usepackage{subfigure}

\theoremstyle{definition}

\newcommand{\drop}[1]{}
\newcommand{\no}{\noindent}
\newcommand{\pt}{{\partial_t}}
\newcommand{\eps}{\varepsilon}

\newcommand{\R}{\mathbb{R}}

\newcommand{\fer}[1]{(\ref{#1})}
\newcommand{\qtext}[1]{\quad\text{#1}}
\newcommand{\abs}[1]{| #1 |}
\newcommand{\nor}[1]{\| #1 \|}
\newcommand{\grad}{\nabla}
\def\O{\Omega}

\newcommand{\B}{{\mathcal{B}}}

\DeclareMathOperator{\Div}{div}

\begin{document}

\maketitle

\begin{abstract}
We use a deterministic particle method to produce numerical approximations to the solutions 
of an evolution cross-diffusion problem for two populations. 

According to the values of the diffusion parameters related to the intra and inter-population repulsion
intensities, the system may be classified in terms of an associated matrix. 
When the matrix is definite positive, the problem is well posed and the Finite Element approximation
produces convergent approximations to the exact solution.  

A particularly important case arises when the matrix is only positive semi-definite and the initial data 
are segregated: the contact inhibition problem. In this case, the solutions may be discontinuous and hence
the (conforming) Finite Element approximation may exhibit instabilities in the neighborhood of the discontinuity.

In this article we deduce the particle method approximation to the general cross-diffusion problem and apply it
to the contact inhibition problem. We then provide some numerical experiments comparing the results produced by 
the Finite Element and the particle method discretizations.
\end{abstract}

\begin{itemize}
 \item[] {\small \emph{Keywords: }
Cross-diffusion system, contact inhibition problem, deterministic particle method, finite element method, numerical simulations. }
 
 \item[] {\small \emph{AMS: }35K55, 35D30, 92D25. }
\end{itemize}

%


\section{Introduction}

In this article, we shall study the numerical approximation to the following problem: Given a fixed $T>0$ and a bounded set $\O\subset\R^m$,  find
$u_i:\O\times(0,T)\to\R$ such that, for $i=1,2$,
\begin{align}
& \pt u_i-\Div J_i(u_1,u_2)=f_i(u_1,u_2) && \qtext{in }Q_T=\O\times(0,T), 
  \label{eq:pde}\\
& J_i(u_1,u_2)\cdot n =0 && \qtext{on }\Gamma_T=\partial\O\times(0,T),
\label{eq:bc}\\
& u_{i}(\cdot,0)=u_{i0} && \qtext{in }\O,	\label{eq:id} 
\end{align}
with flow and competitive Lotka-Volterra functions given by
\begin{align}
 & J_i(u_1,u_2) = u_i\big(a_{i1}\grad u_1 +a_{i2}\grad u_2 +b_i q \big) +
 c_i\grad u_i ,  \label{def:flow}\\
& f_i(u_1,u_2) = u_i \big(\alpha_{i}-\beta_{i1} u_1 - \beta_{i2} u_2 \big) ,
\label{def:reaction}
\end{align}
where, for  $i,j=1,2$, the coefficients $a_{ij},~c_i,~\alpha_i,~\beta_{ij}$ are non-negative constants, $b_i$ is a real constant, $q:\O\times(0,T)\to\R^m$, and $u_{i0}$ are 
non-negative functions.

Problem \fer{eq:pde}-\fer{def:reaction} is a generalization of the cross-diffusion model introduced by Busenberg and Travis \cite{busenberg83} and Gurtin and Pipkin \cite{gurtin84} to take into account the effect
of over-crowding repulsion on the population dynamics, see \cite{gs13} for the modelling
details.

Under the main condition
\begin{equation}\label{H:def_pos}
4 a_{11} a_{22} - (a_{12}+a_{21})^2 > a_0 \qtext{a.e. in } Q_T . 
\end{equation}
for some positive constant $a_0$, it was proven in  \cite{gs13}  the existence of weak solutions
for rather general conditions on the data problem.
Notice that \fer{H:def_pos} implies the following ellipticity condition on the matrix $A=(a_{ij})_{i,j=1}^2$:
\begin{equation*}\label{cond:ellipticity}
\xi^T A~\xi \geq a_0 \nor{\xi}^2  \qtext{for all } \xi\in\R^2.
\end{equation*}
This condition allows us to, through a procedure of approximation, justify the use of
$\ln u_i$ as a test function in the weak formulation of \fer{eq:pde}-\fer{eq:id}.
Then we get, for the entropy functional 
\begin{equation*}
E(t)= \sum_{i=1}^2\int _{\Omega} F(u_i(\cdot,t)) \geq 0, \qtext{with }F(s)=s (\ln s-1)+1,
\end{equation*}
the identity 
\begin{align*}
E(t)+ \int_{Q_t } \Big( \sum_{i=1}^2 ( a_{ii} |\grad u_i|^2 + 2c_i\abs{\grad
\sqrt{u_i}}^2 )+ (a_{12}+a_{21}) \grad u_1\cdot \grad u_2 \Big) \\
  = E(0)+
 \int_{Q_t } \sum_{i=1}^2  \Big( -b_i \, q \cdot \grad u_i +
 f_i(u_1,u_2) \ln u_i \Big) . 
\end{align*}
From  \fer{H:def_pos} and other minor assumptions one then obtains the entropy inequality
\begin{align*}
E(t)+ a_0 \int _{Q_t} ( |\grad u_{ 1}|^2 + |\grad u_{ 2}|^2)\leq (E(0)+C_1) \, \mathrm{e} ^{C_2 t} ,
\end{align*}
providing the key  $L^2(0, T;H^1(\O­))$ estimate of $u_1$ and $u_2$ which allows to prove the existence of weak 
solutions.  

However, it was also proven in \cite{gs13} that condition \fer{H:def_pos} is just a sufficient condition, and that 
solutions may exist for the case of semi-definite positive matrix $A$. 

A particular important case captured by problem \fer{eq:pde}-\fer{def:reaction} is the \emph{contact inhibition problem}, arising in tumor modeling, see for instance Chaplain et al. \cite{chaplain06}. In this case, matrix $A$ is semi-definite positive, and the initial data, describing the spatial distribution of normal and 
tumor tissue,  satisfy 
$\{u_{10}>0\}\cap \{u_{20}>0\} =\emptyset$. 

This free boundary problem was 
mathematically analyzed by Bertsch et al. 
for one \cite{bertsch85} and several spatial dimensions \cite{bertsch12} by using regular Lagrangian flow techniques. In \cite{gs13}, a different approach based on viscosity perturbations was used to prove the existence of solutions. In \cite{gsv13}, the Lagrangian techniques of \cite{bertsch12} were  generalized showing, in particular, the non-uniqueness in the construction of solutions by this method.  

In \cite{gs13}, a conforming Finite Element Method was used both for proving the existence of solutions
of the viscosity approximations (or for the case of $A$ satisfying \fer{H:def_pos}), and for
the numerical simulation of solutions. Since in the case of positive semi-definite   
 matrix $A$ solutions may develope discontinuities in finite time, the FEM approximations
 exhibit instabilities in the neighborhood of these points. In this article we use 
 a deterministic particle method to give an alternative for the numerical simulation of solutions.
 
In our context, deterministic particle methods were introduced, for the scalar linear diffusion equation, by Degond and Mustieles in \cite{degond1989,degond1990}. 
Lions and Mas-Gallic, in \cite{lions01}, gave a rigorous justification of the
method with a generalization to nonlinear diffusion. Gambino et al. \cite{gambino09} studied a particle approximation to a cross-diffusion problem closely related to ours.

Let us finally remark that cross-diffusion parabolic systems have been used 
to model a variety of phenomena since the seminal work of Shigesada, Kawaski 
and Teramoto \cite{shigesada79}. These models
ranges from ecology \cite{gilad07,tian10,gv11,galiano12,sherrat00},
to semiconductor theory \cite{chen07} or granular materials \cite{aranson02,gjv03}, among others.
Global existence and regularity results for the evolution problem 
\cite{kim84,deuring87,yagi93,ggj01,ggj03,chen04,bendahmane10,galiano12},
and for the steady state \cite{lou-ni_96a,lou-ni_96b} have been provided. 
Other interesting properties, such as pattern formation, has been 
studied in \cite{gambino12,ruiz12,gambino13}. Finally, the numerical discretization has received 
much attention, and several schemes have been proposed \cite{ggj01,ggj03,barret04,gambino09,andreianov11,berres11}.


\section{The particle method}

Consider a system of $2N$ particles described by ther masses, $w_k^i>0$, and their 
trajectories, $x_k^i:[0,T]\to\O\subset\R^m$, for $k=1,\ldots,N$ (particle labels) and $i=1,2$ (populations). 
On one hand, the individual  particles state may be modeled by Dirac delta measures
\begin{equation*}
 \epsilon_{x^i_j(t)}(B)=\left\{
\begin{array}{ll}
 1 & \text{if }x^i_k(t)\in B\\
0 &  \text{if }x^i_k(t)\notin B
\end{array}
\right. \qtext{for all } B\in \B(\R^m ),
\end{equation*}
where $\B(\R^m )$ denotes the Borel $\sigma-$algebra generated by open sets in $\R^m$. On the other hand, 
the collective behavior of the discrete system may be given in 
terms of the spatial distribution of particles at time $t$, 
expressed through the \emph{empirical measures}
\begin{equation*}
 u^N_{i}(t)=\sum_{j=1}^{N} w_k^i \epsilon_{x^i_{k}(t)} ,
\end{equation*}
which give the spatial relative frequency of particles of the $i$-th population, at time $t$.


Using the weak formulation of \fer{eq:pde}-\fer{eq:id} and assuming, for the moment, $f_i=0$, we find that $x_k^i$ satisfy, formally
\begin{equation}
\label{ode.1}
 \frac{d}{dt} x_k^i(t)=\frac{1}{u_i^N(x_k^i(t))}J_i(u_1^N(x_k^i(t)),u_2^N(x_k^i(t))).
\end{equation}
We may then initialize this system of ODE's with a suitable particle approximation of the initial data $u_{i0}$ to find an approximation to the solution of problem \fer{eq:pde}-\fer{eq:id}.

Since the Dirac delta measure is difficult to handle, one usually introduces a regularizing non-negative
symmetric kernel $\xi_\eps$ such that $\xi_\eps(x)=\eps^{-m}\xi(x/\eps)$, 
with $\int_{\R^m} \xi =1$ and $\xi_\eps\to \delta$ as $\eps\to 0$, where $\delta$ denotes the Dirac delta distribution. The Gaussian kernel is a common choice.

Then, we look for functions
\begin{equation}
\label{def:ueps}
 u_i^{N,\eps}(x,t)=\sum_{k=1}^N w_k^i\xi_\eps(x-x^i_k(t))
\end{equation}
such that \fer{ode.1} is satisfied with $u_i^N$ replaced by  $u_i^{N,\eps}$, i.e.
\begin{equation}
 \label{ode.2}
 \left\{
 \begin{array}{ll}
 \dfrac{d}{dt} x_k^i(t)= -\Big(a_{i1}\grad u_1^{N,\eps}(x_k^i(t)) +a_{i2}\grad u_2^{N,\eps}(x_k^i(t))+ c_i\dfrac{\grad u_i^{N,\eps}(x_k^i(t))}{u_i^{N,\eps}(x_k^i(t))} +b_i q  \Big), \\
 x_k^i(0)=x_{0k}^i,
 \end{array}
 \right.
 \end{equation}
for $i=1,2$, and with $x_{0k}^i$ such that $u_i^{N,\eps}(x,0)$ is a suitable approximation to the initial data $u_{i0}(x)$ . In the formula above, the gradients are computed as
\begin{equation*}
 \grad u_i^{N,\eps}(x,t)=\sum_{k=1}^N w_k^i\grad \xi_\eps(x-x^i_k(t)).
\end{equation*}

Observe that he existence and uniqueness of solution to problem \fer{ode.2} is guaranteed by the Lipschitz 
continuity of the right hand side term.

\subsection{Implementation}

For the implementation details, we consider the one-dimensional spatial case, i.e. $m=1$. Let us introduce some notation.

Let $\O=(a,b)$. Problem \fer{ode.2} may be rewritten as: 
For $i=1,2$, find $\mathbf{x}^i(t):=(x^i_k(t))_{k=1}^N:[0,T]\to\O^N$ and $\mathbf{w}^i:= (w^i_k)_{k=1}^N$ with $w^i_k>0$,  satisfying 
\begin{equation}
\label{systemODEs}
 \left\{
 \begin{array}{l}
 \dfrac{d\,}{dt} \mathbf{x}^i (t) = \mathbf{g}^i(t,\mathbf{x}^1(t),\mathbf{w}^1,\mathbf{x}^2(t),\mathbf{w}^2) \quad \text{for } t\in (0,T)\, ,\\[1ex]
 \mathbf{x}^i(0) = \mathbf{x}^i_0 ,
 \end{array}
 \right. 
\end{equation}
where 
\begin{equation*}
\begin{array}{ll}
g^i_k(s,\mathbf{y}^1,\mathbf{v}^1,\mathbf{y}^2,\mathbf{v}^2):= & 
-\Big(a_{i1}\grad u^{N,\eps}(y_k^i;\mathbf{y}^1,\mathbf{v}^1) +a_{i2}\grad u^{N,\eps}(y_k^i;\mathbf{y}^2,\mathbf{v}^2)\Big.\\
& \qquad \Big. + c_i\frac{\grad u^{N,\eps}(y_k^i;\mathbf{y}^i, \mathbf{v}^i)}{u^{N,\eps}(y_k^i;\mathbf{y}^i, \mathbf{v}^i)} +b_i q(y^i_k,s)  \Big), 
\end{array}
\end{equation*}
for any $\mathbf{y}^1,\mathbf{v}^1,\mathbf{y}^2,\mathbf{v}^2\in\mathbb{R}^N$ and $s\in (0,T)$, and 
we introduced the following generalized counterpart of (\ref{def:ueps}):
\begin{equation*}
 u^{N,\eps}(x;\mathbf{y},\mathbf{v})=\sum_{k=1}^N v_k \, \xi_\eps(x-y_k) \quad \forall \mathbf{y},\mathbf{v}\in\mathbb{R}^N,\, x\in\mathbb{R}\,.
\end{equation*}
 Notice that, so defined, function $\mathbf{g}^i$ 
 may lead to numerical instabilities when  $u^{N,\eps}(y_k^i;\mathbf{y}^i, \mathbf{v}^i)$ cancels; thus, to avoid divisions by zero, we introduce a parameter $\tilde{\eps}>0$ small enough and approximate $g^i_k(s,\mathbf{y}^1,\mathbf{v}^1,\mathbf{y}^2,\mathbf{v}^2)$ by
\begin{equation*}
\begin{array}{ll}
\tilde{g}^{i}_k(s,\mathbf{y}^1,\mathbf{v}^1,\mathbf{y}^2,\mathbf{v}^2):= & 
-\Big(a_{i1}\grad u^{N,\eps}(y_k^i;\mathbf{y}^1,\mathbf{v}^1) +a_{i2}\grad u^{N,\eps}(y_k^i;\mathbf{y}^2,\mathbf{v}^2)\Big.\\
& \qquad \Big. + c_i\frac{ u^{N,\eps}(y_k^i;\mathbf{y}^i, \mathbf{v}^i) \, \grad u^{N,\eps}(y_k^i;\mathbf{y}^i, \mathbf{v}^i)}{(u^{N,\eps}(y_k^i;\mathbf{y}^i, \mathbf{v}^i)) ^2 + \tilde{\eps}^2} +b_i q(y^i_k,s)  \Big). 
\end{array}
\end{equation*}

In order to approximate system (\ref{systemODEs}), we apply a time discretization 
based on an \emph{implicit midpoint formula}. More precisely, we fix a constant time  step $\Delta t=T/M$ ($M\in\mathbb{N}$) and, given ${u}^{i,n}\approx {u}^i(t_n)$ a particle approximation with the associated positions $\mathbf{x}^{i,n}\approx \mathbf{x}^i(t_n)$ ($i=1,2$) and weights $\mathbf{w}^i $, we approximate $u^{i,n+1}\approx u^i(t_{n+1})$ as follows:
\begin{enumerate}
 \item[1.] Approximate  the position of the associated particles at time $t_{n+1/2}$, $\mathbf{x}^{i,n+1/2}\approx \mathbf{x}^i(t_{n+1/2})$, using  an implicit Euler rule:
\begin{equation}\label{nonlinear-algeq}
\mathbf{x}^{i,n+1/2}=\mathbf{x}^{i,n}+\frac{\Delta t}{2}\, \tilde{\mathbf{g}}^i (t_{n+1/2},\mathbf{x}^{1,n+1/2}, \mathbf{w}^1,\mathbf{x}^{2,n+1/2},\mathbf{w}^2) \, .
\end{equation}
 \item[2.] Approximate  the position of the associated particles at time $t_{n+1}$, $\mathbf{x}^{i,n+1}\approx \mathbf{x}^i(t_{n+1})$, using  an explicit Euler rule: 
$$
\mathbf{x}^{i,n+1}=\mathbf{x}^{i,n+1/2}+\frac{\Delta t}{2} \, \tilde{\mathbf{g}}^i(t_{n+1/2},\mathbf{x}^{1,n+1/2},\mathbf{w}^1,\mathbf{x}^{2,n+1/2},\mathbf{w}^2)\, .
$$
\item[3.] Approximate $u^{i,n+1}\approx u^i(t_{n+1})$ taking advantage of the approximate position of the associated particles $\mathbf{x}^{i,n+1}$:
$$
u^{i,n+1}:=u^{N,\eps}(\cdot ;\mathbf{x}^{i,n+1},\mathbf{w}^i)\, .
$$
\end{enumerate}
Notice that the first step above requires the determination of the solution of the  nonlinear algebraic equation (\ref{nonlinear-algeq}). We compute an approximation to such solution by applying a fixed point algorithm. 
It consists of the following steps:
\begin{enumerate}
 \item[1.1.] Initialize $\mathbf{x}^{i,n+1/2,1}:=\mathbf{x}^{i,n}$.
 \item[1.2.] For $j=1,2,\ldots $, given $\mathbf{x}^{i,n+1/2,j}$,
  $$
\mathbf{x}^{i,n+1/2,j+1}=\mathbf{x}^{i,n}+\frac{\Delta t}{2} \, \tilde{\mathbf{g}}^i(t_{n+1/2},\mathbf{x}^{1,n+1/2,j},\mathbf{w}^1,\mathbf{x}^{2,n+1/2,j},\mathbf{w}^2) \quad (i=1,2) \, .
$$
 \item[1.3.] Check the stopping criteria, 
 $$
 \max_{i=1,2}\max_{x\in C_N } | u_i^{N,\eps}(x;\mathbf{x}^{i,n+1/2,j+1},\mathbf{w}^i) - u_i^{N,\eps}(x;\mathbf{x}^{i,n+1/2,j},\mathbf{w}^i) |    \leq tol \cdot \Delta t \, .
 $$
\end{enumerate}
Here, $tol>0$ is a tolerance parameter and $C_N$ is a finite set of points. In our experiments, we take $tol=4\cdot 10^{-6}$, and $C_N$ given by a uniform grid of points of $\Omega$, 
$$
C_N :=\{x_k:=a+(k-1)\Delta x;\quad k=1,2, \ldots , N \} \, ,
$$
with $\Delta x=\frac{b-a}{N-1}$.

Notice that the initial condition in (\ref{systemODEs}) is given in terms of the initial location of the particles, whereas that of the original problem, (\ref{eq:id}), provides a value of their spatial distribution. In this sense, given an initial condition $u^i(\cdot ,0)=u_{i0}$ in $\O$, we must study how to initialize $\mathbf{x}^i_0$ and $\mathbf{w}^i$ in such a way that
\begin{equation}\label{icond_aprox}
 u^{N,\eps}(\cdot ;\mathbf{x}^i_0,\mathbf{w}^i)\approx u_{i0}\quad\text{in }\O\, .
\end{equation}
On the one hand, to initialize the  positions, we used a uniform grid of $\O$, that is, we simply took  the set of points $C_N$:
$$
\mathbf{x}^i_0:=(x_k)_{k=1}^N \, .
$$
On the other hand, to initialize the weights, one is temted to impose (\ref{icond_aprox}) exactly on the points of $C_N$. However, in doing so, negative weights arise, spoiling the convergence of the method. In consequence, we compared two different strategies:
\begin{itemize}
\item define $w^i_k:=\Delta x\, u_{i0}(x_k)$;
\item solve the following constrained linear least-squares problem:
$$
\text{find $\mathbf{w}^i\in\mathbb{R}_+^N$ that minimizes $||A\mathbf{W}-\mathbf{u}_{i0}||_2^2$ for $\mathbf{W}\in\mathbb{R}_+^N$ }\, ,
$$
where $A_{k,l}:= \xi_\eps(x_k-x_l)$ and $\mathbf{u}_{i0}:= (u_{i0}(x_k))_{k=1}^N$.
\end{itemize}
Although the first approach is faster, it introduces too much diffussion, whereas the second provides a more accurate approximation while preserving stability. 

A similar strategy applies to the particle redistribution  after several time steps. Indeed, let us recall that particle redistribution is tipically needed in particle simulations to avoid that particles get concentrated in some parts of the domain: if so, too big gaps between particles will arise in other parts, producing numerical instabilities. We refer to \cite{barbaetal2005} for a review on this issue and some alternative strategies that apply.

Finally, the boundary conditions (\ref{eq:bc}) are taken into account by means of a specular reflection whenever a particle location is changed. Although this is the most common approach to non-flow boundary conditions, 
border effects are noticiable in simple discretization schemes.


\section{Numerical experiments}

\subsection{The Finite Element Method discretization}

For comparison purpouses, we use the Finite Element approximation of 
\fer{eq:pde}-\fer{eq:id} introduced in \cite{gs13}, eventually using the regularization 
of the flows \fer{def:flow} given by 
\begin{equation}
\label{reg:flow}
J_i^{(\delta)}(u_1,u_2)= J_i(u_1,u_2) +\frac{\delta}{2} \Delta(u_iu),
\end{equation}
for $\delta\ge0$. The approximation is obtained by using a 
semi-implicit Euler scheme in time and a $\mathbb{P}_1$ continuous finite element
approximation in space, see
\cite{gs13} for the details where, in particular, the convergence of the fully discrete approximation to the 
continuous solution is proved.

We now sketch the FEM scheme. Let   $\Delta t>0$ be the time step of the discretization. For $t=t_0=0$, set $u_{\epsilon
i}^0=u_i^0$. Then, for $n\geq 1$ the problem is to find $u_{\epsilon
i}^{n}\in S_h$, the finite element space of piecewise $\mathbb{P}_1$-elements, such that for,  $i=1,2$,
\begin{equation}\label{eq:pde_discr.s4}
\begin{array}{l}
\frac{1}{\Delta t}\big( u^n_{\epsilon i}-u^{n-1}_{\epsilon i} , \chi )^h + \big(
J^{(\delta)}_i(\Lambda _{\epsilon } (u^n_{\epsilon 1}),\Lambda _{\epsilon } (u^n_{\epsilon
2}),\nabla u^n_{\epsilon 1},\nabla u^n_{\epsilon 2}  ) ,\nabla\chi \big)^h =\\ [2ex]
\hspace*{1cm} = \big(\alpha_{i} u^n_{\epsilon i} - \lambda _{\epsilon } (u^n_{\epsilon i}) (
\beta_{i1} \lambda _{\epsilon } (u^{n-1}_{\epsilon 1}) + \beta_{i2} \lambda
_{\epsilon } (u^{n-1}_{\epsilon 2}) ) , \chi \big)^h , 
\end{array}
\end{equation}
for every $ \chi\in S^h$.
Here, $(\cdot,\cdot)^h$ stands for a discrete semi-inner product on
$\mathcal{C}(\overline{\Omega} )$. The parameter $\epsilon>0$ makes reference to the
regularization introduced by some functions $\lambda_\epsilon$ and $\Lambda_\epsilon$, which
converge to the identity as $\epsilon\to0$, see \cite{gs13}.

Since \eqref{eq:pde_discr.s4} is a nonlinear algebraic problem, we use a fixed point
argument to approximate its solution,  $(u_{\epsilon 1}^n,u_{\epsilon 2}^n)$, at each time
slice $t=t_n$, from the previous approximation $u_{\epsilon i}^{n-1}$.  Let $u_{\epsilon
i}^{n,0}=u_{\epsilon i}^{n-1}$. Then, for $k\geq 1$ the problem is to find $u_{\epsilon
i}^{n,k}$ such that for $i=1,2$, and for all $\chi \in S^h$
\begin{equation*}
\begin{array}{l}
 \frac{1}{\Delta t}\big( u^{n,k}_{\epsilon i}-u^{n-1}_{\epsilon i} , \chi )^h
+ \big( J^{(\delta)}_i(\Lambda _{\epsilon } (u^{n,k-1}_{\epsilon 1}),\Lambda _{\epsilon }
(u^{n,k-1}_{\epsilon 2}),\nabla u^{n,k}_{\epsilon 1},\nabla u^{n,k}_{\epsilon 2}  )
,\nabla\chi \big)^h =\\ [2ex] \hspace*{1cm} = \big(\alpha_{i} u^{n,k}_{\epsilon i} - \lambda
_{\epsilon } (u^{n,k-1}_{\epsilon i}) ( \beta_{i1} \lambda _{\epsilon } (u^{n-1}_{\epsilon
1}) + \beta_{i2} \lambda
_{\epsilon } (u^{n-1}_{\epsilon 2}) ) , \chi \big)^h . 
\end{array}
\end{equation*}
We use the stopping criteria $\max _{i=1,2}
\|u_{\epsilon,i}^{n,k}-u_{\epsilon,i}^{n,k-1}\|_\infty <\text{tol}$, for empirically chosen
values of $\text{tol}$, and set $u_i^n=u_i^{n,k}$.

%
%


\subsection{Experiments}

\no\emph{Experiment 1.} 
We consider a particular situation of the contact-inhibition problem in which an explicit
solution of \eqref{eq:pde}  may be computed in terms of a suitable combination of the
Barenblatt explicit solution of the porous medium equation, the Heavyside function and the
trajectory of the contact-inhibition point. To be precise, we construct a solution to the
problem
\begin{align}
& \partial_t u_{i} - (u_i (u_1+u_2)_x)_x =0&& \text{in }(-R,R)\times(0,T)=Q_T,&\label{eq:s1}\\
&u_i (u_1+u_2)_x = 0 &&\text{on } \{-R,R\}\times(0,T),& \label{eq:s2}
\end{align}
with
\begin{equation}
\label{def:us}
 u_{10}(x)=H(x-x_0)B(x,0), \quad u_{20}(x)=H(x_0-x)B(x,0).
\end{equation}
Here, $H$ is the Heavyside function and $B$ is the Barenblatt solution of the porous medium
equation corresponding to the initial datum $B(x,-t^*)=\delta_0$, i.e.

\begin{equation}
\notag \label{barenblatt}
 B(x,t)= 2 (t+t^*)^ {-1/3} \big[1-\frac{1}{12}x^ 2(t+t^*)^ {-2/3}\big]_+.
\end{equation}
For simplicity, we consider problem \eqref{eq:s1}-\eqref{def:us} for $T>0$ such that
$r(T)<R^ 2$, with $r(t)=\sqrt{12}(t+t^*)^ {1/3}$, so that $B(R,t)=0$ for all $t\in[0,T]$.
The point $x_0$ is the initial contact-inhibition point, for which we assume $|x_0|<r(0)$,
i.e. it belongs to the interior of the support of $B(\cdot,0)$, implying that the initial
mass of both populations is positive. Observe that $u_1+u_2$ satisfies the porous medium 
equation, implying regularity properties for this sum, among others, the differentiability 
in the interior of its support. 

It can be shown \cite{gsv13} that 
the functions
\begin{equation}
\notag
 u_1(x,t)=H(x-\eta(t))B(x,t),\quad u_2(x,t)=H(\eta(t)-x)B(x,t),
\end{equation}
with $\eta(t)=x_0(t/t^*)^{1/3}$, are a weak solution of problem \eqref{eq:s1}-\eqref{def:us}.

We use the FEM scheme and the particle method scheme to produce approximate solutions to 
problem \fer{eq:s1}-\fer{def:us} for $R=1$ and a resolution of $N=1000$ (nodes or particles). In this experiment, 
the FEM scheme behaves well without the addition of the regularizing term in \fer{reg:flow}, i.e. we
take $\delta=0$. The initial data
is given by \fer{def:us} with $t_*=0.01$. 

We run the experiments till the final time $T=0.01$ is reached. We use a small
time resolution in order to capture the discontinuity of the exact solution. As suggested in \cite{degond1990}
the restriction $\Delta t = C\eps^2$ must be impossed in order to get stability. We chose
\begin{equation*}
 \Delta t = 0.1 \eps^2,\qtext{with }\eps=0.15 dx^{0.75},
\end{equation*}
and $dx=1/N$. This high time resolution implies that the fixed point algorithms to solve the nonlinearities
is scarcely used. Particle spatial redistribution is neither needed in this experiment.

\begin{figure}[ht]
\centering
 \subfigure%
 {\includegraphics[width=6.25cm,height=4.75cm]{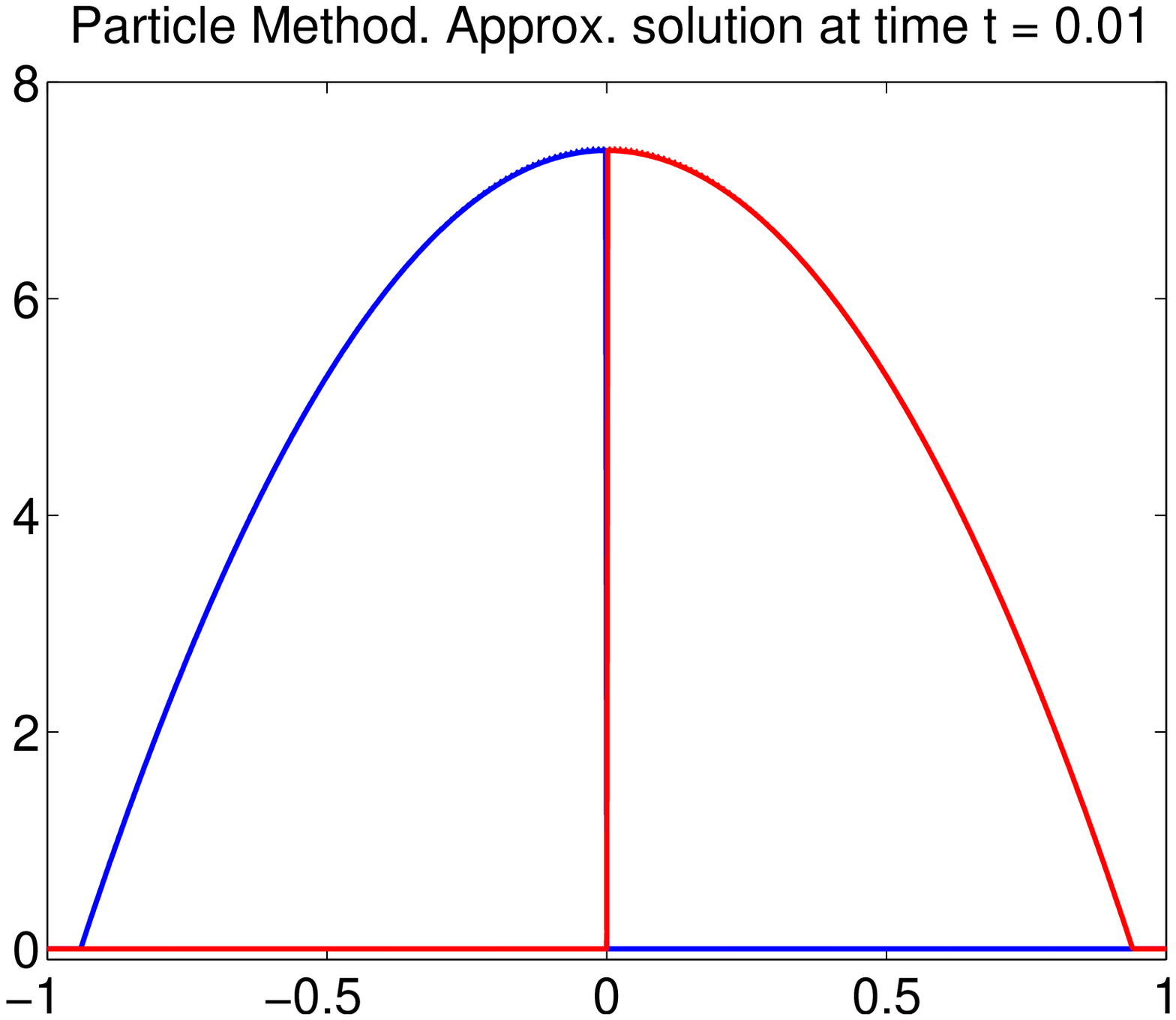}}
 \subfigure
 {\includegraphics[width=6.25cm,height=4.75cm]{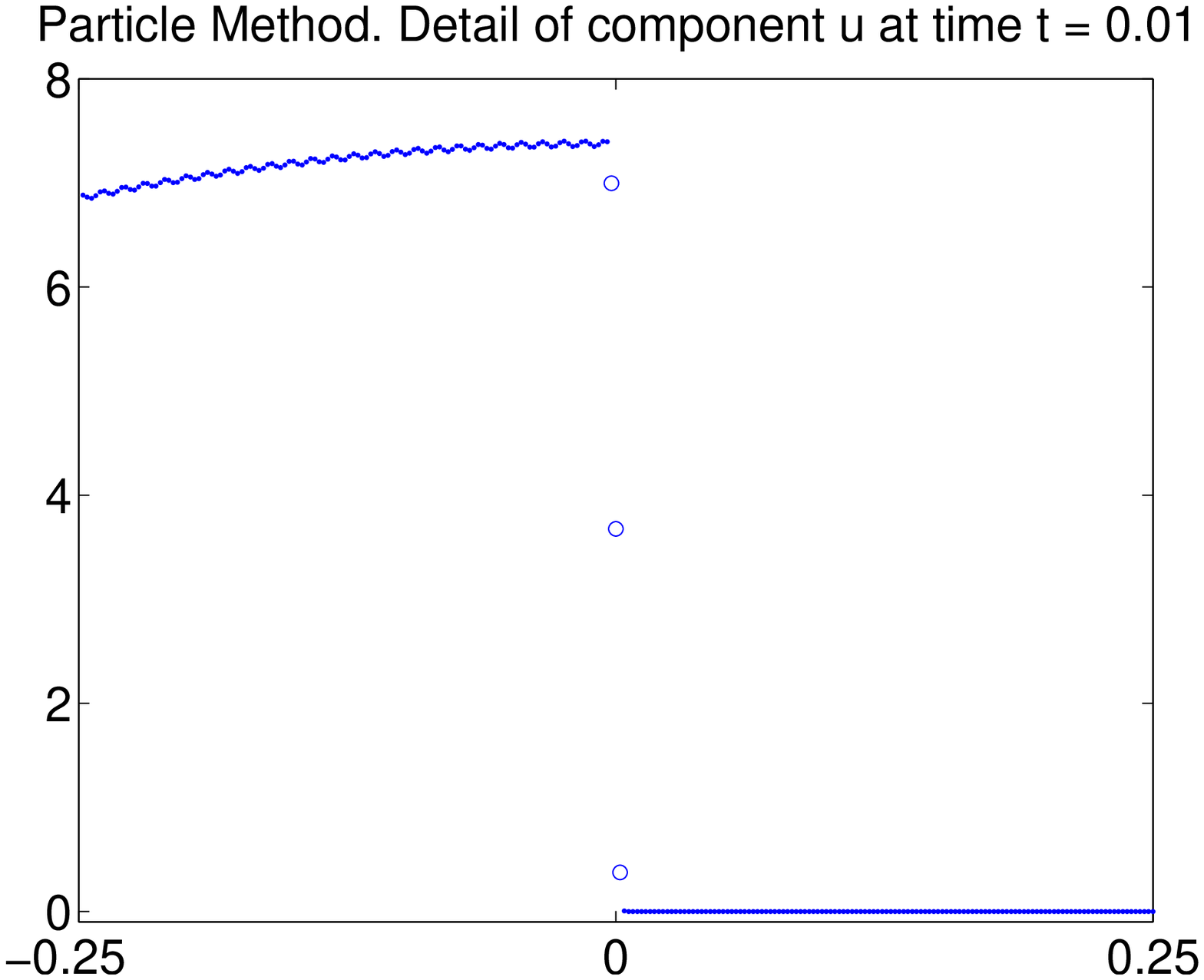}} 
 \caption{{\small Experiment 1. Left panel: Exact and particle method approximation. At this scale, both are indistinguishable. Right panel: A detail of the solution showing the oscillating instabilities and the need 
 of only three particles to covering the discontinuity of the exact solution.}} 
\label{exp1_1.fig}
\end{figure}

\begin{figure}[ht]
\centering
 \subfigure%
 {\includegraphics[width=6.25cm,height=4.75cm]{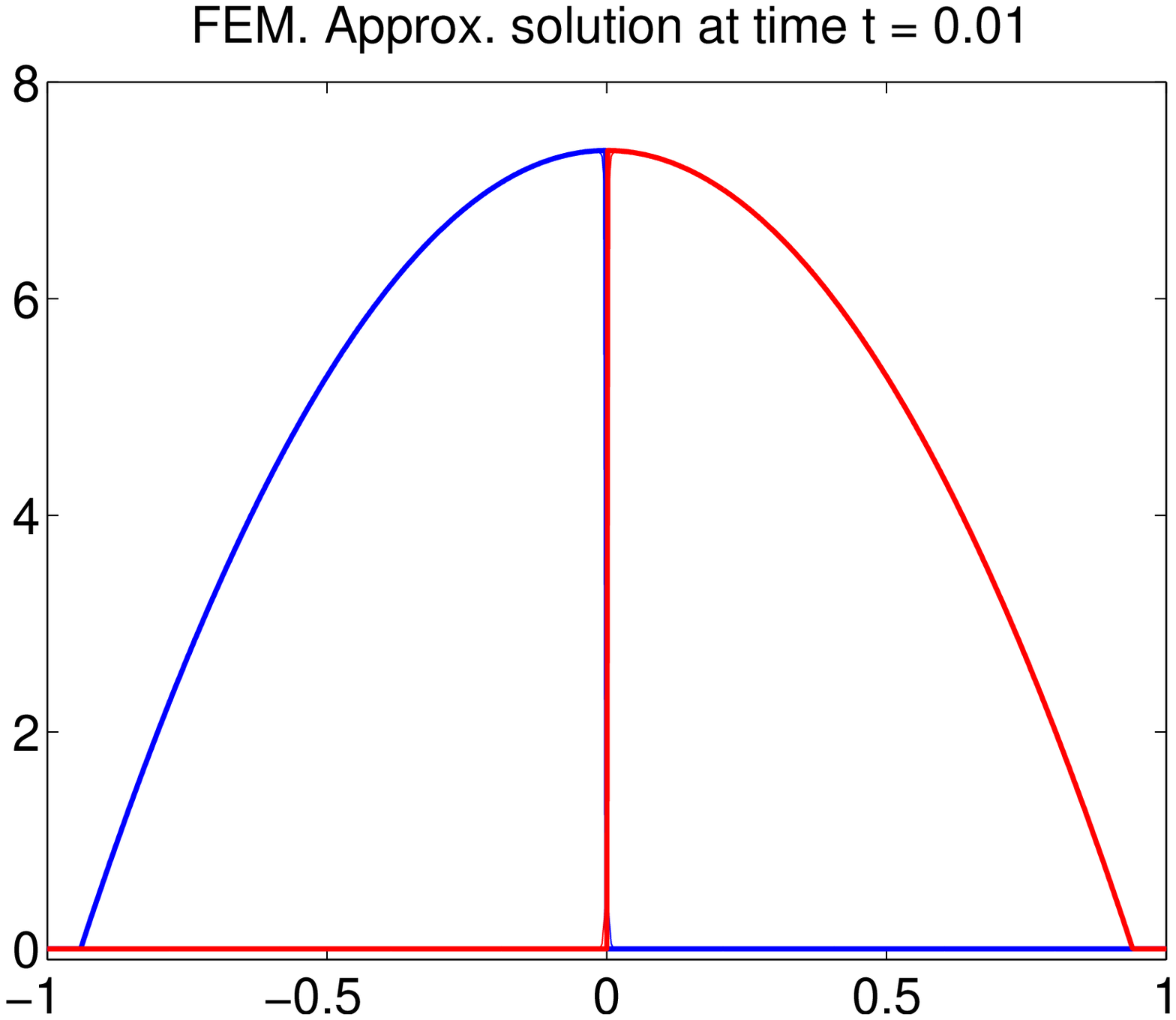}}
 \subfigure
 {\includegraphics[width=6.25cm,height=4.75cm]{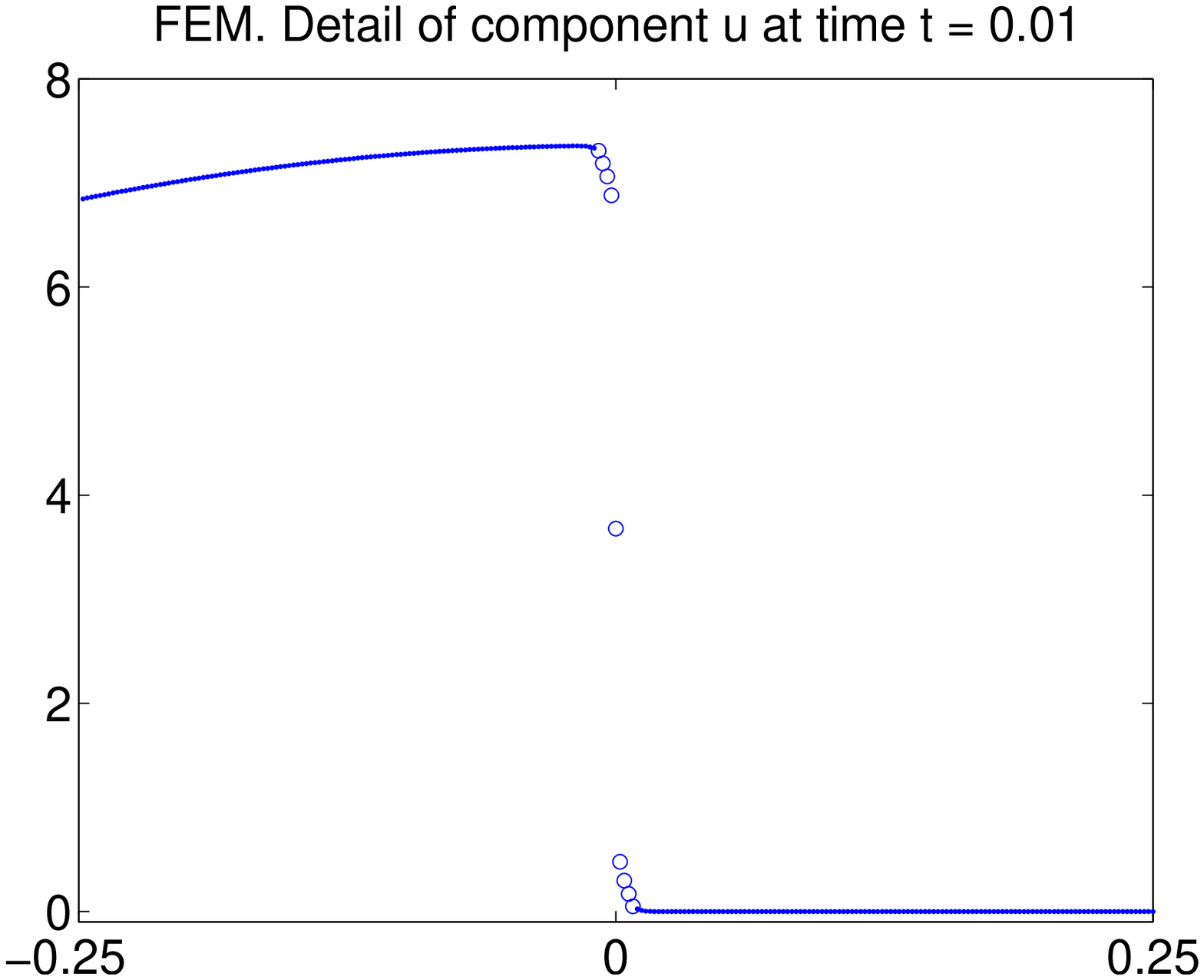}} 
 \caption{{\small Experiment 1. Left panel: Exact and FEM approximation. Again, at this scale, both are indistinguishable. Right panel: A detail of the solution showing that more nodes than in the particle method are needed to cover the discontinuity of the exact solution. However, no instabilities are present.}} 
\label{exp1_2.fig}
\end{figure}

Although both algorithms produce similar results, i.e. a good approximation outside a small 
neighborhood of the discontinuity $x=0$, they behave in a different way. On one hand, the particle method needs
 fewer particles to cover the discontinuity, see Figs.~\ref{exp1_1.fig} and \ref{exp1_2.fig}. On the other hand, the particle method creates oscillating instabilities in a large region of the positive part of the solution, effect which is not observed in the case of the FEM. In any case, the global errors are similar. In particular, the mean relative square error is of  order $10^{-2}$.

\bigskip 

\begin{figure}[ht]
\centering
 \subfigure%
 {\includegraphics[width=6.25cm,height=4.75cm]{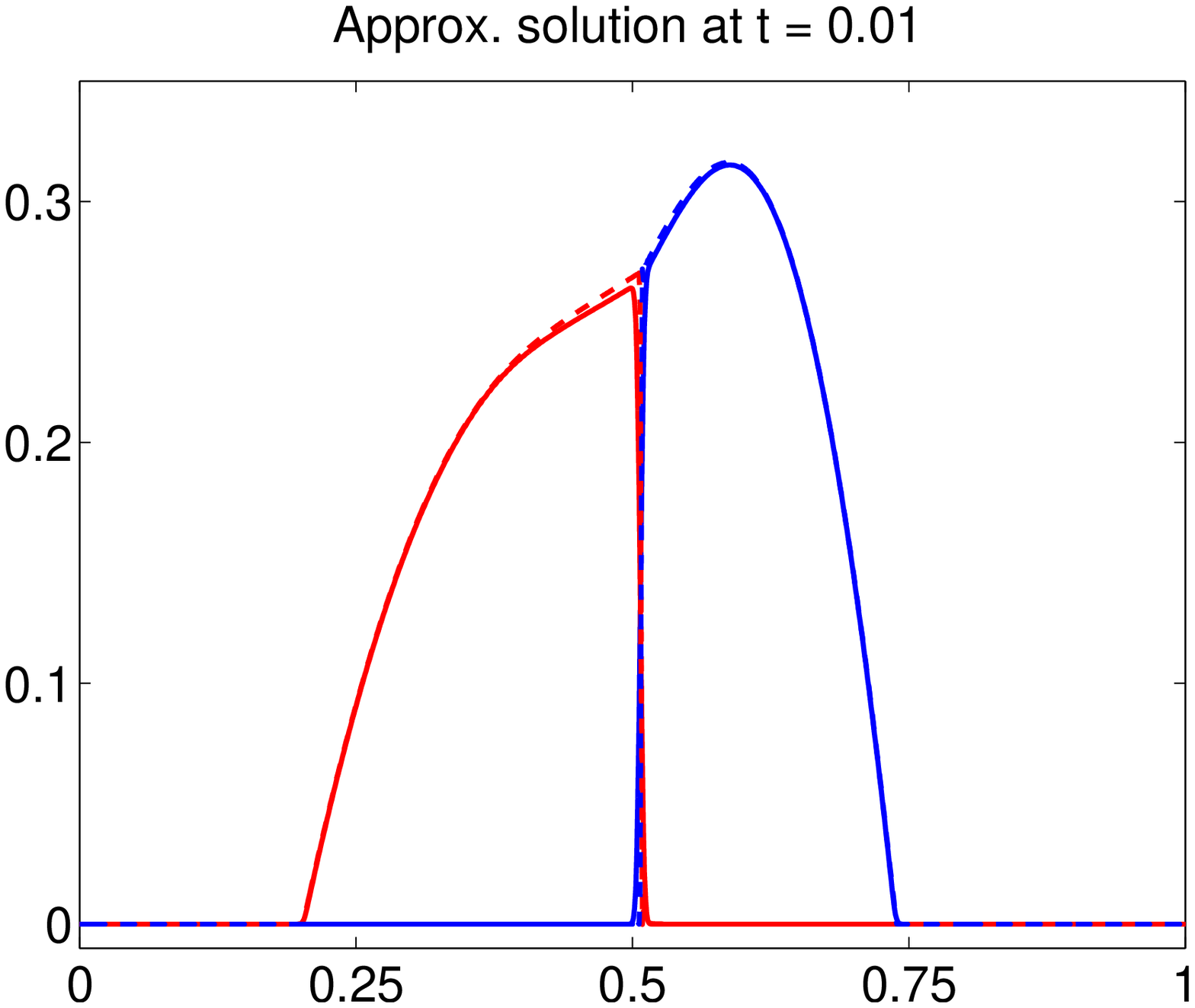}}
 \subfigure
 {\includegraphics[width=6.25cm,height=4.75cm]{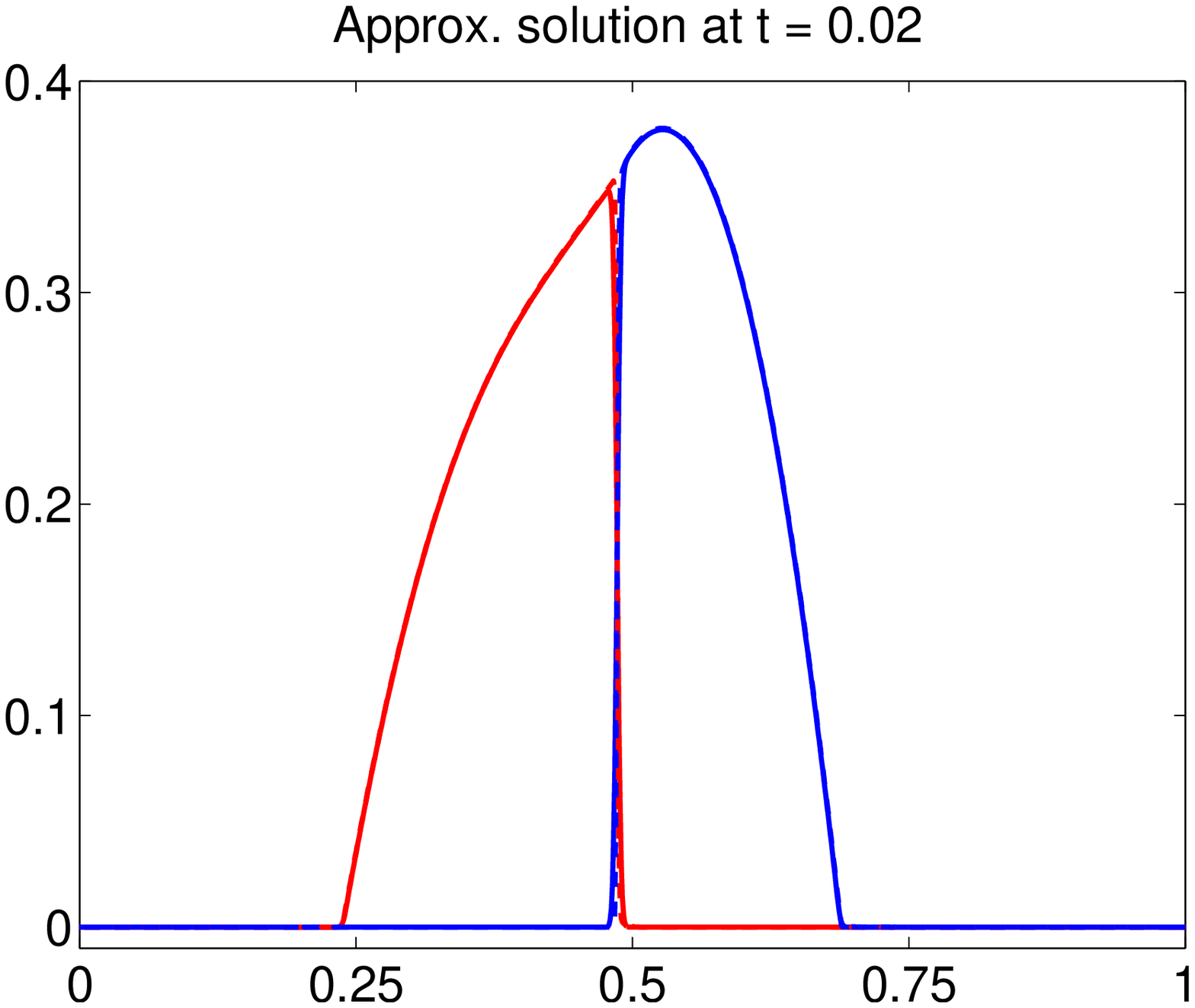}} 
 \caption{{\small Experiment 2. Particle (continuous line) and FEM (dotted line) approximations. Only at the discontinuity point a difference between both approximations may be observed. }} 
\label{exp2_1.fig}
\end{figure}

\no\emph{Experiment 2.} 
Another instances of the contact-inhibition problem are investigated. In Fig.~\ref{exp2_1.fig}, we show approximate 
transient solutions obtained by the Particle method (continuous line) and Finite Element method (dotted line)
of two problems given in the form
\begin{align*}
& \partial_t u_{i} - (u_i (a_{i1}u_{1x}+a_{i2}u_{2x} +b_iq) )_x =0&& \text{in }(0,1)\times(0,T),&\\
&u_i (a_{i1}u_{1x}+a_{i2}u_{2x} +b_iq) = 0 &&\text{on } \{0,1\}\times(0,T),&
\end{align*}
with
\begin{equation*}
 u_{i0}(x)=\text{e}^{ - (\frac{x-x_i}{0.001})^2},\qtext{for } x_1=0.4,\quad x_2=0.6,
\end{equation*}
and $q(x)=-3(x-0.5)$.

For the first problem (left panel of Fig.~\ref{exp2_1.fig}) we choose 
\begin{equation*}
 A=(a_{ij})=\left(\begin{array}{cc} 3 & 3 \\ 1 & 1\end{array}\right),\quad b_1=b_2=0,
\end{equation*}
and for the second problem (right panel of Fig.~\ref{exp2_1.fig}) we choose 
\begin{equation*}
 A=(a_{ij})=\left(\begin{array}{cc} 1 & 1 \\ 1 & 1\end{array}\right),\quad b_1=1,~ b_2=10,
\end{equation*}
implying that matrix $A$ is positive semi-definite in both cases. Differently than in Experiment 1, 
the sum $u_1+u_2$  does not satisfy more than a continuity regularity property for these problems. Indeed, a jump of the derivative 
may be observed at the contact inhibition point. This effect is explained by the differences in the flows on
the left and on the right of the contact inhibition point. 

As it can be seen in the figures, Particle and Finite Element methods provide a similar approximation.
Only at the contact inhibition point some differences may be observed. The $L^2$ relative error is, as in the Experiment~1,  of  order $10^{-2}$.

The discretization parameters 
are 
\begin{equation*}
 \Delta t = 0.5 \eps^2,\qtext{with }\eps=0.5 dx^{0.75},
\end{equation*}
and $dx=1/N$, with $N=1000$. The regularized flow \fer{reg:flow} is used in the FEM scheme with $\delta=10^{-3}$. The experiments are run till the final time $T=0.01$ (left panel) and $T=0.02$ (right panel). 

\section{Conclusions}

The contact inhibition problem, i.e. Equations \fer{eq:pde}-\fer{def:reaction} completed with 
initial data with disjoint supports, is an interesting problem from the mathematical point of 
view due mainly to the possibility of their solutions developing discontinuities in finite time.

Although there has been some recent progress in the analytical understanding of the problem \cite{bertsch12,gs13,gsv13}, the numerical analysis is still an open problem. In this paper we have presented
in some extension a Particle method to produce approximations to the solutions of the problem. We have 
compared these solutions to exact solutions and to approximate solutions built through the Finite Element Method.

In general terms, the Particle Method is more computer time demanding, and somehow unstable with respect to the resolution parameters when we compare to FEM. However, if the resolution is high enough, it can better capture the discontinuities arising in the solution. Although we performed our experiments in a one dimensional setting, 
Particle Methods are specially useful for higher dimensions, due to the easiness of their implementation and parallelization. In future work we shall investigate these extensions.

\newpage

\begin{figure}[ht]
\centering
 \subfigure%
 {\includegraphics[width=6.25cm,height=4.75cm]{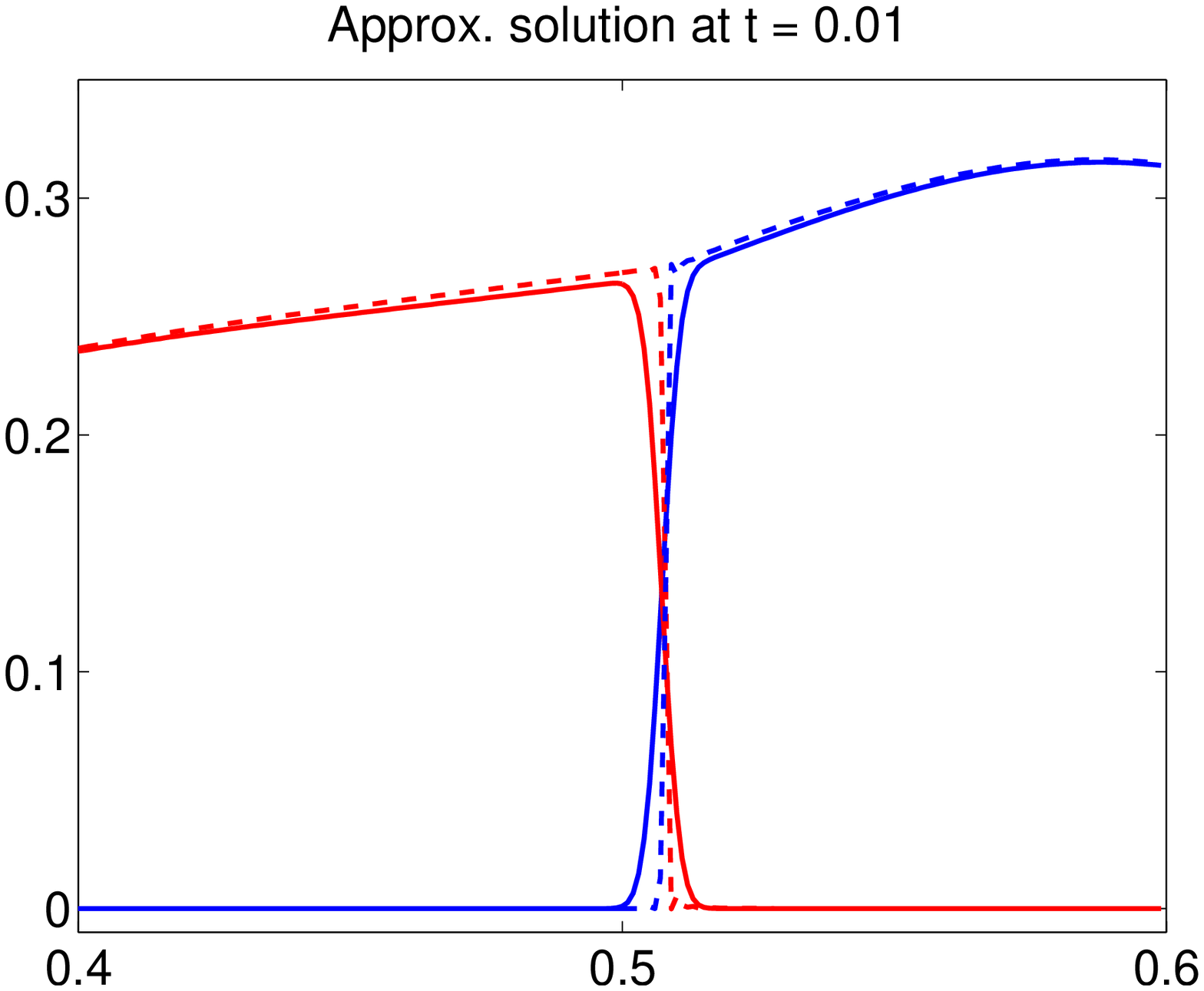}}
 \subfigure
 {\includegraphics[width=6.25cm,height=4.75cm]{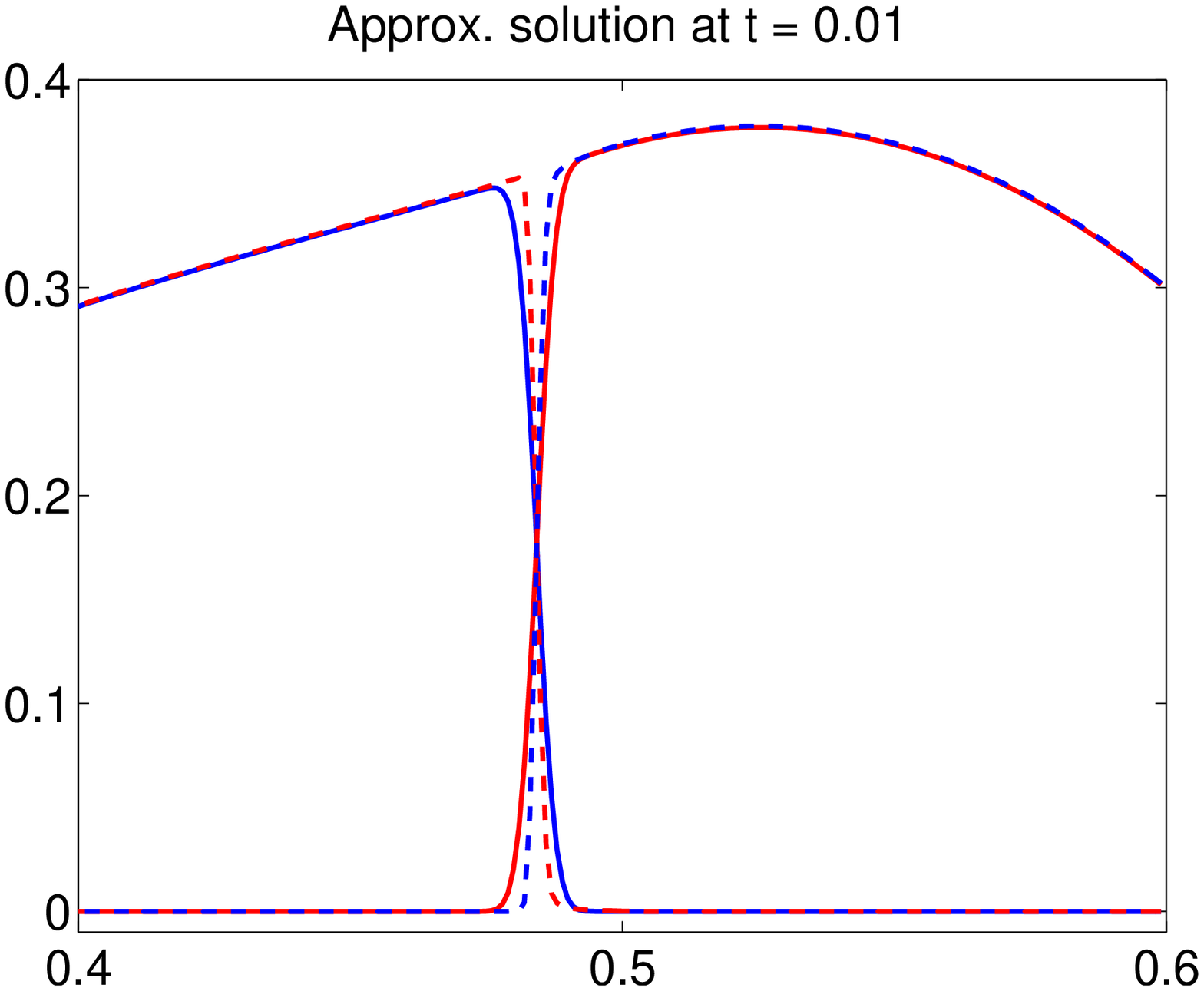}} 
 \caption{{\small Experiment 2. Detail of Fig.~\ref{exp2_1.fig} at the contact inhibition point. }} 
\label{exp2_2.fig}
\end{figure}


\begin{thebibliography}{10} 


\bibitem{andreianov11}
M. Andreianov, B. Bendahmane, R. Ruiz-Baier,
Analysis of a finite volume method for a cross-diffusion model in population dynamics,
Math. Mod. Meth. Appl. Sci. 21(2) (2011) 307--344.

\bibitem{aranson02}
I. S. Aranson, L. S. Tsimring,
Continuum theory of partially fluidized granular flows,
Phys. Rev. E (3) 65 (2002) 061303.

\bibitem{barbaetal2005}
L. A. Barba, A. Leonard, C. B. Allen, Advances in viscous vortex methods-meshless spatial
adaption based on radial basis function interpolation, Int. J. Numer. Meth. Fluids 47 (2005) 387--421.


\bibitem{barret04}{ J.~W. Barrett, J.~F. Blowey}, { Finite element
approximation of a
nonlinear cross-diffusion population model}, Numer. Math. 98 (2004)
 195--221.

 
\bibitem{bendahmane10} { M. Bendahmane}, { Weak and classical solutions to
predator-prey system with cross-diffusion}, Nonlinear Anal. 73 (2010)
 2489--2503.

\bibitem{berres11}
 S. Berres, R. Ruiz-Baier,
 A fully adaptive numerical approximation for a two-dimensional epidemic 
 model with nonlinear cross-diffusion,
 Nonlinear Anal. Real World Appl. 12 (2011) 2888--2903.

 
\bibitem{bertsch85}
{ M. Bertsch, M.~E. Gurtin, D. Hilhorst, L.~A. Peletier}, { On
interacting populations that disperse to avoid crowding: preservation of
segregation}, J. Math. Biol. 23 (1985)  1--13


\bibitem{bertsch12}
{ M. Bertsch, D. Hilhorst, H. Izuhara, M. Mimura}, { A nonlinear
parabolic-hyperbolic system for contact inhibition of cell-growth}, Diff. Equ.
Appl. 4 (2012)  137--157.

\bibitem{busenberg83}
{ S.~N. Busenberg, C.~C. Travis}, { Epidemic models with spatial spread
due to population migration}, J. Math. Biol. 16 (1983)  181--198.

\bibitem{chaplain06}
M. Chaplain, L. Graziano, L. Preziosi, 
Mathematical modelling of the loss of tissue compression responsiveness and its role 
in solid tumour development,
Math. Med. Biol. 23 (2006) 197--229.


\bibitem{chen04}{ L. Chen, A. J\"ungel}, { Analysis of a multidimensional parabolic 
population model with strong cross-diffusion}, SIAM J. Math. Anal.  36(1)  (2004) 301--322.

\bibitem{chen07}{ L. Chen, A. J\"ungel},
Analysis of a parabolic cross-diffusion semiconductor model with electron-hole scattering, 
Comm. Partial Differential Equations, 32 (2007) 127--148.


%
%
%
%
%


\bibitem{degond1989} P. Degond, S. Mas-Gallic, The weighted particle method for convection–diffusion equations. Part I. The case of an isotropic viscosity and
Part II. The anisotropic case, Math. Comp. 53 (1989) 485--525.

\bibitem{degond1990} P. Degond, F.-J. Mustieles, A deterministic approximation of diffusion equations using particles, SIAM J. Sci. Stat. Comput. 11 (1990) 293--
310.

\bibitem{deuring87} { P. Deuring}, { An initial-boundary value problem for a certain
density-dependent diffusion system}, Math. Z. 194 (1987)  375-396.


\bibitem{galiano12}{ G. Galiano}, { On a cross-diffusion population model deduced from
mutation and splitting of a single species}, Comput. Math. Appl. 64(6) (2012) 1927-1936.

\bibitem{ggj01}{ G. Galiano, M.~L. Garz\'on, A. J\"ungel}, { Analysis and
numerical
solution of a nonlinear cross-diffusion system arising in population dynamics}, 
RACSAM Rev. R. Acad. Cienc. Exactas Fís. Nat. Ser. A Mat.  95(2)  (2001)
 281--295.

\bibitem{ggj03}{ G. Galiano, M.~L. Garz\'on, A. J\"ungel},
{ Semi-discretization in
time  and numerical convergence of solutions of a nonlinear cross-diffusion
population model},  Numer. Math.  93(4)  (2003)  655--673. 


\bibitem{gjv03}
{ G. Galiano, A. J\"ungel, J. Velasco},
A parabolic cross-diffusion system for granular materials,
SIAM J. Math. Anal. 35(3) (2003) 561--578.

\bibitem{gs13}
G. Galiano, V. Selgas, 
On a cross-diffusion segregation problem arising from a model of interacting particles.
To appear in Nonlinear Anal. Real World Appl.

\bibitem{gsv13}
G. Galiano, S. Shmarev , J. Velasco,
Existence and non-uniqueness of segregated solutions to a class of cross-diffusion systems. In preparation.




\bibitem{gv11}
G. Galiano, J. Velasco, 
Competing through altering the environment: A cross--diffusion population 
model coupled to transport-darcy flow equations,
Nonlinear Anal. Real World Appl. 12(5) (2011) 2826--2838.


\bibitem{gambino09}
{ G. Gambino, M.C. Lombardo and M. Sammartino,} 
A velocity-diffusion method for a Lotka-Volterra system with nonlinear cross and self-diffusion, 
Appl. Numer. Math. 59 (2009)  1059--1074.

\bibitem{gambino12}
{ G. Gambino, M.C. Lombardo and M. Sammartino,} 
Turing instability and traveling fronts for a nonlinear reaction–diffusion system with 
cross-diffusion Original, Math. Comput. Simul.82(6) (2012) 1112--1132.

\bibitem{gambino13}
{ G. Gambino, M.C. Lombardo and M. Sammartino,} 
Pattern formation driven by cross--diffusion in a 2D domain,
Nonlinear Anal. Real World Appl. 14(3) (2013) 1755--1779.

\bibitem{gilad07}
E. Gilad, J. von Hardenberg, A. Provenzale, M. Shachak, E. Meron.
A mathematical model of plants as ecosystem engineers,
J. Theoret. Biol. 244(4) (2007) 680--691.


\bibitem{gurtin84}
{ M.~E. Gurtin, A.~C. Pipkin}, { On interacting populations that disperse
to avoid crowding}, Q. Appl. Math. 42 (1984)  87--94.

%
%

\bibitem{kim84}{ J.~U. Kim}, { Smooth solutions to a quasi-linear system
of diffusion equations for a certain population model}, Nonlinear Anal. 8 (1984)
 1121--1144.



\bibitem{lions01} P.-L. Lions, S. Mas-Gallic, Une m\'ethode particulaire d\'eterministe pour des équations diffusive non lin\'eaires, C. R. Acad. Sci. Paris S\'er. I
Math. 332 (2001) 369--376.

\bibitem{lou-ni_96a}
 { Y. Lou, W.~M. Ni, Y. Wu}, 
 Diffusion, self-diffusion and cross-diffusion, 
 J. Differ. Equations 131(1) (1996) 79--131.

 
\bibitem{lou-ni_96b} { Y. Lou, W.~M. Ni, Y. Wu}, { The global existence of solutions for
a cross-diffusion system}, Adv. Math. Beijing 25 (1996)  283-284.





\bibitem{ruiz12}
R. Ruiz-Baier, C. Tian,
Mathematical analysis and numerical simulation of pattern formation under 
cross-diffusion,
Nonlinear Anal. Real World Appl. 14(1) (2013) 601--612.

\bibitem{sherrat00}
J. A. Sherratt, 
Wavefront propagation in a competition equation with a new
motility term modelling contact inhibition between cell populations,
R. Soc. Lond. Proc. Ser. A Math. Phys. Eng. Sci.  456 (2000) 2365--2386.


\bibitem{shigesada79}
{ N. Shigesada, K. Kawasaki, E. Teramoto}, { Spatial segregation of
interacting species}, J. Theor. Biol. 79 (1979)  83--99.



\bibitem{tian10}
C. Tian, Z. Lin, M. Pedersen,
Instability induced by cross-diffusion in reaction-diffusion systems,
Nonlin. Anal. RWA 11 (2010) 1036--1045.


%
%
%


\bibitem{yagi93}{ A. Yagi}, { Global solution to some quasilinear
parabolic system in population dynamics}, Nonlinear Anal. 21 (1993)  603--630.

\end{thebibliography}
\end{document}